\documentclass[preprint,3p,authoryear]{elsarticle}

\usepackage[utf8]{inputenc} 
\usepackage[T1]{fontenc}    
\usepackage{amssymb}
\usepackage{amsmath}
\usepackage{hyperref}       
\usepackage{url}            
\usepackage{booktabs}       
\usepackage{amsfonts}       
\usepackage{nicefrac}       
\usepackage{microtype}      
\usepackage{graphicx}
\usepackage{bm}
\usepackage{listings}       
\usepackage{xcolor}
\usepackage{cleveref} 
\graphicspath{ {./images/} }

\hypersetup{ colorlinks=false, linkcolor=black, filecolor=black, urlcolor=black }

\definecolor{codegray}{gray}{0.4}
\definecolor{codekw}{rgb}{0.13,0.29,0.53}
\definecolor{codestr}{rgb}{0.0,0.45,0.2}
\lstdefinestyle{pystyle}{
  language=Python,
  basicstyle=\ttfamily\footnotesize,
  keywordstyle=\color{codekw}\bfseries,
  commentstyle=\color{codegray}\itshape,
  stringstyle=\color{codestr},
  showstringspaces=false,
  columns=fullflexible,
  breaklines=true,
  frame=single,
  framesep=4pt,
  rulecolor=\color{codegray},
  numbers=none,
  aboveskip=6pt,
  belowskip=4pt
}

\newcommand{\Fc}{\bm{F}_c}
\newcommand{\Fv}{\bm{F}_v}
\newcommand{\uvec}{\bm{u}}
\newcommand{\nvec}{\bm{n}}
\newcommand{\grad}{\nabla}
\newcommand{\dd}{\,\mathrm{d}}

\journal{Journal of Open Source Software}

\begin{document}

\begin{frontmatter}


\title{JAX-FVM: A differentiable, entropy-stable finite volume solver on
unstructured meshes for compressible flows}
\tnotetext[label1]{}
\author[INRIA]{Guillaume de Romémont\corref{cor1}}
\ead{guillaume dot romemont at protonmail dot com}
\cortext[cor1]{Corresponding author.}

\affiliation[INRIA]{organization={INRIA},
           addressline={200 Av. de la Vieille Tour},
           city={Talence},
           postcode={33400 },
           country={France}}


\begin{abstract}
We present JAX-FVM, an open-source, fully differentiable finite volume method
(FVM) for the two-dimensional compressible Euler and Navier-Stokes equations on
unstructured triangular meshes. The solver is written entirely in JAX, so that
every operation : mesh connectivity, flux evaluation, slope limiting, and time
integration is just-in-time compiled, vectorised, and end-to-end
differentiable through automatic differentiation (AD), and runs transparently on
CPU or GPU. On the numerical side, JAX-FVM is built around an
entropy-conservative Tadmor/Ismail-Roe two-point flux supplemented with
entropy-variable Rusanov or Roe dissipation, second-order MUSCL reconstruction of
primitive variables with least-squares gradients and Venkatakrishnan limiting,
and a family of explicit (RK2-4) and matrix-free implicit (Newton, SDIRK2)
time integrators whose Jacobian actions are obtained by AD. The combination of an
\emph{unstructured}-mesh compressible FVM with end-to-end differentiability fills
a gap left by existing differentiable CFD frameworks, which are almost
exclusively restricted to structured grids or spectral discretisations. We
describe the governing equations, the discretisation, the software architecture,
and a set of standard verification cases. 
The code is openly available at \url{https://github.com/guigzair/jax_fvm}.

\end{abstract}

\begin{keyword}
differentiable programming \sep JAX \sep finite volume method \sep unstructured
meshes \sep entropy stability \sep compressible Euler and Navier-Stokes equations
\end{keyword}

\end{frontmatter}

\section{Introduction}
\label{sec:intro}
Computational fluid dynamics (CFD) has become an indispensable tool across
science and engineering, from aerodynamics and combustion to biomedical flows
\citep{brunton2020ml}. The canonical use of a CFD solver is a forward one: given
a geometry, boundary conditions, and physical parameters, one integrates the
governing equations to obtain the flow field. Yet many of the most consequential
problems are inverse in nature, since boundary conditions are only partially
observed, constitutive or model parameters are poorly characterised, and
turbulence closures carry model-form uncertainty. Such problems, for instance
recovering vascular resistance and compliance from sparse clinical measurements
\citep{pant2017inverse}, inferring geometry from surface pressure in shape
optimisation, or identifying turbulence-closure corrections from reference data
\citep{duraisamy2019turbulence}, are naturally cast as PDE-constrained
optimisation, and solving them at scale requires the gradient of the solver
output with respect to its inputs. Finite-difference or derivative-free
strategies become prohibitive in high dimension, so scalable sensitivity
computation, historically via adjoint methods, is essential.

Two developments have reshaped how these gradients are obtained. First,
differentiable programming and automatic differentiation (AD) frameworks such as
JAX \citep{jax2018github}, PyTorch \citep{paszke2019pytorch}, TensorFlow
\citep{abadi2016tensorflow}, and Julia \citep{bezanson2017julia} allow the
derivative of an entire program to be evaluated to machine precision without
hand-derived adjoints \citep{sapienza2024diffprog}. Traditional adjoint CFD
implementations, e.g. SU2 \citep{economon2016su2} or dolfin-adjoint
\citep{mitusch2019dolfinadjoint}, remain code-intrusive and are architecturally
decoupled from the AD ecosystems in which modern machine-learning components are
built. Second, the hybridisation of machine learning (ML) with CFD has grown into
a dominant paradigm \citep{brunton2020ml,chen2022physics}. Early hybrid models
embedded offline-trained closures into conventional solvers
\citep{duraisamy2019turbulence}, but such loosely coupled models are not
solver-consistent and can degrade or destabilise once deployed in the live solver
\citep{wu2019illconditioned}. This has motivated fully coupled, differentiable
approaches in which ML components and discrete operators are optimised end-to-end
through the solver, including solver-in-the-loop and learned turbulence closures
\citep{list2022learned,um2020solver}, data-driven discretisations
\citep{barsinai2019learning, de2025data, de2026data}, and hybrid architectures that combine differentiable
PDE solvers with graph neural networks
\citep{belbuteperes2020combining,pfaff2020meshgraphnets}. In all of these,
differentiability through the solver is a structural prerequisite rather than a
convenience.

These requirements have driven a wave of GPU-native differentiable CFD platforms
built on AD frameworks. Representative examples include JAX-CFD
\citep{kochkov2021ml}, a JAX-native incompressible solver widely used in ML-CFD
studies; $\Phi$Flow \citep{holl2020phiflow}, a backend-agnostic differentiable PDE
framework; JAX-Fluids \citep{bezgin2023jaxfluids}, a high-order finite-volume
solver for compressible and two-phase flows; WaterLily
\citep{weymouth2024waterlily}, a Julia immersed-boundary solver; PICT
\citep{franz2025pict}, a differentiable multi-block PISO solver; and Diff-FlowFSI
\citep{fan2026diffflowfsi,fan2024fsi} for turbulence and fluid-structure
interaction. A common limitation, however, is that these frameworks operate almost
exclusively on \emph{structured} Cartesian grids, where the regular data layout
maps cleanly onto vectorised array operations. Unstructured meshes, which are
indispensable for the body-conforming discretisation of geometrically complex
domains and are the workhorse of industrial CFD, remain largely absent from the
differentiable-simulation ecosystem, because their irregular connectivity and
variable neighbour counts are hostile to the SIMD/SIMT execution model of
accelerators. A recent preprint, DiFVM \citep{du2026difvm}, addresses this gap for
incompressible flow by recasting the finite-volume operators as graph
message-passing primitives \citep{pfaff2020meshgraphnets}; to our knowledge,
however, no code is publicly available, and the compressible, entropy-stable
regime is not its focus.

On the numerical side, robustly simulating compressible flows on unstructured
meshes raises its own difficulties. Nonlinear conservation laws admit
discontinuous (shock) solutions, and naive discretisations are prone to spurious
oscillations and nonlinear instability. Two complementary ingredients address
this. High-resolution finite-volume schemes recover second-order accuracy away
from discontinuities through gradient reconstruction and slope limiting on
unstructured stencils \citep{barth1989design,venkatakrishnan1993convergence},
while entropy-stable schemes provide a nonlinear stability framework by enforcing
a discrete entropy inequality. Building on Tadmor's entropy-conservative flux
theory \citep{tadmor1987entropy,tadmor2003entropy} and affordable
entropy-consistent Euler fluxes \citep{ismail2009affordable}, together with the
classical Roe approximate Riemann solver \citep{roe1981approximate}, one obtains
schemes that are provably entropy-stable and robust for strong shocks such as the
forward-facing step and double Mach reflection \citep{woodward1984numerical} and
the two-dimensional Riemann configurations of \citet{lax1998solution}. Casting the
entire scheme in an AD framework is particularly attractive here, since the
entropy-variable Jacobians underlying the dissipation and the implicit solvers can
be formed exactly by differentiation rather than by hand.

JAX-FVM sits at the intersection of these two threads. It is, to our knowledge,
one of the first openly available differentiable finite-volume solvers for the
compressible Euler and Navier-Stokes equations that operates natively on
unstructured triangular meshes, coupled with an entropy-stable Tadmor/Ismail-Roe
discretisation. Relative to the concurrent DiFVM effort \citep{du2026difvm}, the
present work is complementary: it targets compressible flows with entropy-stable
numerics in a compact, research-oriented implementation. 
\section{Statement of need}
\label{sec:need}
JAX-FVM is designed to serve two communities that are rarely addressed by a
single code. To the CFD practitioner it offers a compact, entropy-stable
compressible solver on unstructured meshes; to the scientific-machine-learning
practitioner it offers an unstructured compressible solver that is
\emph{differentiable end to end}, so that gradients of any scalar functional of
the flow can be propagated back to the initial condition, boundary data, physical
parameters, or an embedded neural model.

Three properties, taken together, distinguish JAX-FVM from existing tools:
\begin{enumerate}
  \item \textbf{Differentiability.} The entire residual--connectivity gathers,
  reconstruction, flux, and boundary treatment--is written in JAX and is
  differentiable by construction. AD is used not only to expose gradients to the
  user, but internally, to form Jacobian-vector products for the implicit time
  integrators and the entropy-consistent dissipation operators.
  \item \textbf{Entropy stability.} The default flux is an entropy-conservative
  Tadmor/Ismail-Roe two-point flux \citep{tadmor1987entropy,ismail2009affordable}
  augmented with entropy-variable dissipation, so that the semi-discrete scheme
  respects a discrete entropy inequality--a property that is difficult to obtain
  and verify without an AD-enabled implementation.
  \item \textbf{Unstructured meshes.} Existing differentiable CFD frameworks such
  as JAX-CFD \citep{kochkov2021ml}, JAX-Fluids \citep{bezgin2023jaxfluids}, and
  $\Phi$Flow \citep{holl2020phiflow} operate on Cartesian/structured grids. To
  our knowledge there is little to no publicly available differentiable
  \emph{unstructured} FVM code for compressible flow. JAX-FVM targets exactly
  this niche.
\end{enumerate}

The solver is intended as a research platform.

\section{Governing equations}
\label{sec:equations}
JAX-FVM solves the two-dimensional compressible Navier-Stokes equations in
conservative form,
\begin{equation}
  \partial_t \mathbf{w} + \grad \cdot \big(\Fc(\mathbf{w}) - \Fv(\mathbf{w}, \grad \mathbf{w})\big) = 0,
  \qquad
  \mathbf{w} = (\rho,\ \rho u,\ \rho v,\ E)^\top ,
  \label{eq:ns}
\end{equation}
where $\rho$ is the density, $\uvec = (u,v)$ the velocity, $p$ the pressure, and $E$ the total
energy per unit volume. The inviscid (convective) flux is
\begin{equation}
  \Fc(\mathbf{w}) =
  \begin{pmatrix}
    \rho \uvec \\
    \rho \uvec \otimes \uvec + p\,\mathbb{I} \\
    (E + p)\,\uvec
  \end{pmatrix},
\end{equation}
closed by the ideal-gas relation. In the implementation the pressure is recovered
from a non-dimensional form that carries a reference Mach number $M$,
\begin{equation}
  E = \frac{p}{(\gamma-1)\,M^2} + \dfrac{1}{2}\rho\,(u^2+v^2),
  \label{eq:eos}
\end{equation}
with heat-capacity ratio $\gamma = 1.4$ by default. The viscous flux uses a Newtonian 
stress tensor under the Stokes hypothesis (zero bulk viscosity) together with Fourier heat conduction,
\begin{align}
  \bm{\tau} &= \mu\Big(\grad\uvec + \grad\uvec^\top
              - \tfrac{2}{3}(\grad\cdot\uvec)\,\mathbb{I}\Big),
  &
  \bm{q} &= -k\,\grad T,
  &
  T &= \frac{p}{\rho R},
  \label{eq:viscous}
\end{align}
so that $\Fv \cdot \nvec = (0,\ \bm{\tau}\cdot \nvec,\ (\bm{\tau}\uvec + \bm{q})\cdot\cdot \nvec)^\top$,
with dynamic viscosity $\mu$, thermal conductivity $k$, and specific gas constant
$R$. The compressible Euler equations are obtained as the inviscid limit
$\Fv \equiv 0$.

\section{Finite volume discretisation on unstructured meshes}
\label{sec:fvm}

\subsection{Mesh representation}
The computational domain is triangulated with the MeshPy \footnote{\url{https://github.com/inducer/meshpy}}/Triangle Delaunay mesh
generator \citep{shewchuk1996triangle}, subject to a maximum-area and a
minimum-angle ($30^\circ$ by default) constraint. The \texttt{Mesh} class stores
the vertices, the triangle-vertex table, cell barycentres and areas, the global
face list with per-face lengths and outward normals, and the cell-to-cell
adjacency (\texttt{neighbors}) together with a face-connectivity map. Periodic
boundaries are handled by matching opposite boundary faces and rewiring the
adjacency accordingly. Each boundary face carries an integer marker that selects
its boundary condition (Section~\ref{sec:bc}); mesh output is written to VTK
through \texttt{meshio}.

\subsection{Semi-discrete form}
Integrating \eqref{eq:ns} over a triangle $C_i$ and applying the divergence
theorem yields the cell-centred semi-discrete update
\begin{equation}
  \frac{\dd \mathbf{w}_i}{\dd t} = -\,\mathcal{R}_i(\mathbf{w}),
  \qquad
  \mathcal{R}_i(\mathbf{w}) = \frac{1}{|C_i|} \sum_{f \in \partial C_i}
  \hat{\bm{H}}\big(\mathbf{w}_i^f, \mathbf{w}_{j(f)}^f; \nvec_f\big)\, \ell_f,
  \label{eq:semidiscrete}
\end{equation}
where $\mathbf{w}_i$ is the cell-averaged state, $|C_i|$ the cell area, $\ell_f$ the face
length, $\nvec_f$ the outward unit normal, and $\hat{\bm{H}}$ a consistent
numerical flux evaluated from the reconstructed left/right face states
$\mathbf{w}_i^f, \mathbf{w}_{j(f)}^f$. The residual $\mathcal{R}_i$ is the quantity that is
just-in-time compiled and differentiated throughout the code.

\subsection{Gradient reconstruction and slope limiting}
Second-order accuracy is obtained by a MUSCL reconstruction of the
\emph{primitive} variables $\mathbf{u} = (\rho,u,v,p)$. Cell gradients are computed by
a weighted least-squares (LSQ) fit that minimises
\begin{equation}
  \sum_{j \in \mathcal{N}(i)} w_{ij}\,
  \big(\grad \mathbf{u}_i \cdot \Delta\mathbf{x}_{ij} - (\mathbf{u}_j - \mathbf{u}_i)\big)^2 ,
  \qquad w_{ij} = |\Delta\mathbf{x}_{ij}|^{-2},
  \label{eq:lsq}
\end{equation}
where $\Delta\mathbf{x}_{ij}$ joins the barycentres of cells $i$ and $j$; the
resulting $2\times2$ normal equations are solved per cell and per variable with a
doubly-vectorised (\texttt{vmap}) linear solve. The face states are then
extrapolated as
\begin{equation}
  \mathbf{u}_i^f = \mathbf{u}_i + \phi_i\,\big(\grad\mathbf{u}_i \cdot (\mathbf{x}_f - \mathbf{x}_i)\big),
\end{equation}
with a slope limiter $\phi_i \in [0,1]$. Two limiters are provided: the smooth
Venkatakrishnan limiter \citep{venkatakrishnan1993convergence},
\begin{equation}
  \phi(a,b) = \frac{a^2 + 2ab + \omega}{a^2 + 2b^2 + ab + \omega},
  \qquad \omega = (K\,h)^3,
  \label{eq:venka}
\end{equation}
(with $h$ a local length scale and $K$ a tunable constant) and the classical
\texttt{minmod} limiter. This LSQ+MUSCL construction follows the standard
upwind unstructured meshes methodology of \citet{barth1989design}.

For the viscous flux, the LSQ gradients are evaluated at faces by
\begin{equation*}
    \nabla \mathbf{w}_{ij} = \overline{\nabla \mathbf{w}}_{ij} + \left( \frac{\mathbf{w}_j - \mathbf{w}_i}{|r_j - r_i|} - \overline{\nabla \mathbf{w}}_{ij} \cdot \hat{r}_{ij} \right) \hat{r}_{ij}
\end{equation*}
with 
\begin{equation*}
    \overline{\nabla \mathbf{w}}_{ij} = \dfrac{1}{2}(\nabla \mathbf{w}_i + \nabla \mathbf{w}_j), \quad \hat{r}_{ij} = \dfrac{r_j - r_i}{|r_j - r_i|}
\end{equation*} 
Even if the stencil needed to compute the gradient is larger than the one needed for the convective flux, the scheme avoids the odd-even decoupling problem.

\subsection{Numerical fluxes}
\label{sec:flux}
The scheme is formulated in terms of the entropy variables associated with the
specific entropy $s = \ln(p/\rho^\gamma)$. The default convective flux is the
entropy-conservative two-point flux of Tadmor
\citep{tadmor1987entropy,tadmor2003entropy} in the affordable Ismail-Roe form
\citep{ismail2009affordable}. Writing the Ismail-Roe parameter vector
$\bm{z} = \sqrt{\rho/p}\,(1,\,u,\,v,\,p)^\top$ and denoting by
$\overline{(\cdot)}$ the arithmetic and by $(\cdot)^{\ln}$ the logarithmic mean
(evaluated with a numerically stable Taylor expansion near equal arguments), the
consistent, entropy-conservative interface state $(\hat\rho,\hat u,\hat v,\hat p)$
is reconstructed and assembled into the physical flux $\hat{\bm{H}}^{\mathrm{EC}}$.
A dissipation term is then added,
\begin{equation}
  \hat{\bm{H}} = \hat{\bm{H}}^{\mathrm{EC}}
  - \tfrac{1}{2}\,\alpha\,|\lambda_{\max}|\,
    \frac{\partial \mathbf{w}}{\partial \bm{\eta}}\,(\bm{\eta}_R - \bm{\eta}_L),
  \label{eq:dissipation}
\end{equation}
where $\bm{\eta}$ are the entropy variables and
$\lambda_{\max} = c/M + |\uvec\cdot\nvec|$ is the maximum wave speed. The jump is
mapped from entropy to conservative variables through the Jacobian
$\partial \mathbf{w}/\partial\bm{\eta}$, whose action is evaluated \emph{exactly} by a
forward-mode AD (\texttt{jax.jvp}) call rather than an analytical
hand-derivation; this yields an entropy-stable Rusanov (local Lax-Friedrichs)
dissipation. The dimensionless coefficient $\alpha$ tunes the amount of numerical
viscosity. As alternatives, a full Roe matrix-dissipation flux
\citep{roe1981approximate} with Roe-averaged eigenstructure and
entropy-variable scaling of the wave strengths and a plain central/Rusanov flux
are also available.

For the Navier-Stokes equations the viscous flux \eqref{eq:viscous} is evaluated
at faces from face-averaged primitive states and face-averaged LSQ gradients, and
subtracted from the convective contribution in the residual.

\section{Boundary conditions}
\label{sec:bc}
Boundary conditions are encoded through the integer face markers attached to the
mesh, and are applied by constructing an appropriate ghost/exterior face state:
marker~1 denotes a periodic face; marker~2 a solid wall; marker~3 a supersonic
(fully prescribed) inlet; marker~4 an outlet (zero-gradient extrapolation); and
marker~5 a subsonic inlet. Walls are treated either as inviscid slip walls
(velocity mirrored about the face normal) or as no-slip walls, selected at run
time.

\section{Time integration}
\label{sec:time}
Given the semi-discrete system \eqref{eq:semidiscrete}, several time integrators
are provided. The default explicit scheme is the two-stage strong-stability
Heun/RK2 method; classical RK3 or RK4 are also available. For stiff or steady problems,
two implicit schemes are implemented: a backward-Euler solve and a singly
diagonally-implicit RK2 (SDIRK2, with coefficient $x = 1 - 1/\sqrt{2}$). Each
implicit stage solves the nonlinear system
\begin{equation}
  \bm{G}(\mathbf{w}) \equiv \mathbf{w} - \mathbf{w}^n + \Delta t\,\mathcal{R}(\mathbf{w}) = 0
\end{equation}
by a few Newton iterations. Crucially, the required Jacobian-vector products are
never assembled: the action
$\bm{v} \mapsto \bm{v} + \Delta t\,(\partial\mathcal{R}/\partial\mathbf{w})\,\bm{v}$ is
computed matrix-free by forward-mode AD (\texttt{jax.jvp}) and passed to a GMRES
solver (\texttt{jax.scipy.sparse.linalg.gmres}), with light under-relaxation for
robustness. The explicit time-step size is set from a convective CFL condition
based on the local spectral radius; for viscous runs a diffusive limit
$\Delta t \sim \Delta x^2/\nu$ is also enforced, and the minimum of the two is
taken.

\section{Differentiability and implementation}
\label{sec:impl}
JAX-FVM is written entirely against \texttt{jax.numpy}. The residual
\eqref{eq:semidiscrete} and every time-stepper are wrapped in \texttt{jax.jit},
with the (static) mesh object passed as a compile-time argument. Per-cell
operations such as the LSQ solve \eqref{eq:lsq} and the face-connectivity search
are expressed as nested \texttt{jax.vmap} maps. Forward-mode AD (\texttt{jax.jvp})
plays a dual role: it supplies the entropy-variable jumps in the dissipation
\eqref{eq:dissipation} and the Jacobian actions in the implicit solvers.

Because the whole residual is differentiable, reverse-mode AD can propagate the
gradient of any scalar objective; for instance a mismatch against reference data
or an aerodynamic functional back to the initial condition, the boundary data,
the physical parameters $(\mu,k,\gamma,\dots)$, or the parameters of an embedded
model. The library depends on \texttt{flax} and \texttt{optax}, so that such
optimisation and inverse-design loops can be assembled directly on top of the
solver. Execution is hardware-agnostic: the same code runs on CPU or GPU (a CUDA
build of JAX is used for GPU acceleration), and double precision can be enabled
when strict conservation checks are required. A suite of diagnostics: total
entropy, kinetic energy, enstrophy, palinstrophy, vorticity, and Mach
number supports quantitative verification of the structure-preserving
properties of the scheme.

\section{Software structure and workflow}
\label{sec:workflow}
The package is organised into a mesh module (\texttt{src/mesh}), a solver module
with sub-packages for the compressible Euler, compressible Navier-Stokes,
projection-based incompressible Navier-Stokes, and a pseudo-spectral vorticity
solver (\texttt{src/solvers}), and a library of initial conditions and mesh
geometries (\texttt{src/Cases}). A typical simulation follows four steps: build a
mesh, prescribe an initial condition as primitive variables and convert to
conserved variables, march in time, and post-process. Listing~\ref{lst:api}
condenses the forward-facing-step example.

\begin{lstlisting}[style=pystyle, caption={Minimal driver for the Mach~3
forward-facing step (condensed from \texttt{examples/Euler\_example.ipynb}).},
label={lst:api}]
import jax
from jax_fvm.src.mesh import Mesh_cases
from jax_fvm.src.Cases import Test_Cases
from jax_fvm.src.solvers.Euler import Euler
from jax_fvm.src.solvers import helper

# 1. Build the unstructured mesh (geometry + boundary markers)
mesh = Mesh_cases.Forward_Step().build(h=5e-4)

# 2. Initial condition: primitives (rho, u, v, p) -> conserved state W
Primitives = Test_Cases.ForwardFacingStep().build(mesh)
W = helper.getConserved(Primitives)

# 3. Time integration: stable dt from CFL, then explicit RK2 march
t_final, CFL = 4.0, 0.2
dt = helper.get_dt(W, mesh, CFL=CFL)
N_t = int(t_final / dt) + 1
for n in range(N_t):
    W = jax.jit(Euler.time_step_RK2, static_argnums=(1,))(W, mesh, dt)

# 4. Post-process: back to primitives and plot the fields
Primitives = helper.getPrimitive(W)
mesh.plot_solution(Primitives[..., 0], labels=r'$\rho$')
\end{lstlisting}

\section{Example and verification}
\label{sec:examples}
\begin{figure}[htbp]
  \centering
  \includegraphics[width=0.9\linewidth]{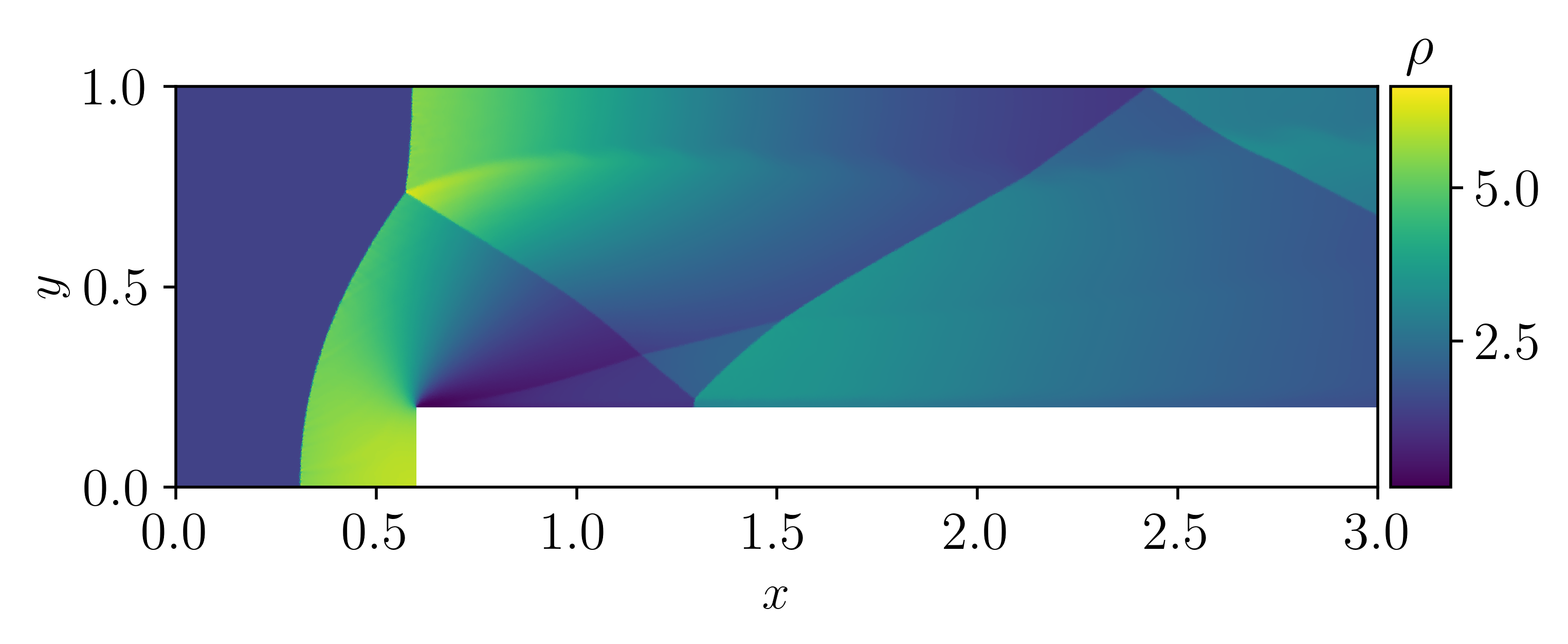}\\[4pt]
  \includegraphics[width=0.9\linewidth]{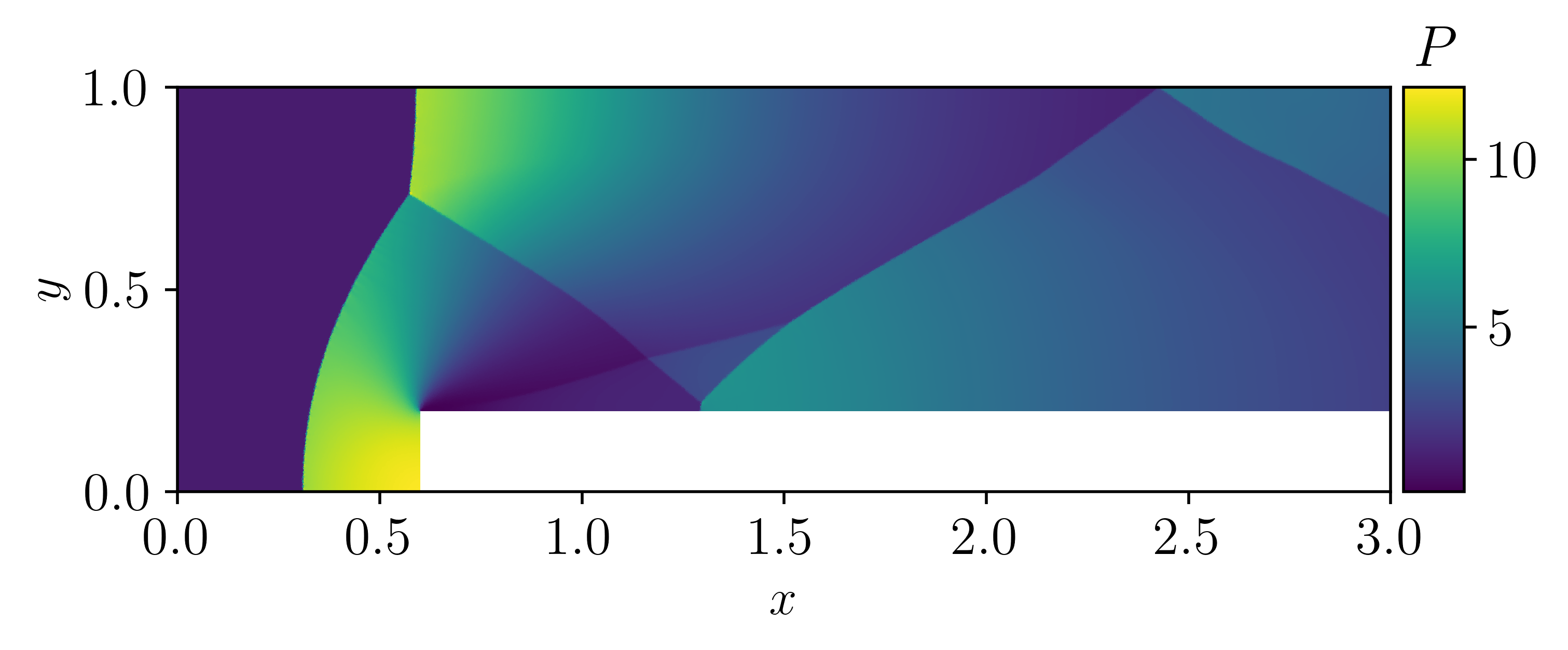}
  \caption{Mach~3 forward-facing step at $t=4$ on an unstructured triangular mesh:
  density $\rho$ (top) and pressure $p$ (bottom). The bow shock, Mach stem,
  reflected shocks, and slip line are resolved. Computed with the
  entropy-conservative Tadmor/Ismail-Roe flux, MUSCL reconstruction with
  Venkatakrishnan limiting, and explicit RK2 time stepping.}
  \label{fig:ffs}
\end{figure}

The headline example is the Mach~3 forward-facing step of
\citet{woodward1984numerical}: a uniform flow at $(\rho,u,v,p) = (1.4,3,0,1)$
enters a wind tunnel containing a forward-facing step. The simulation develops the
canonical unsteady shock system, a detached bow shock ahead of the step, a Mach
stem, reflected shocks, and a contact/slip line which JAX-FVM captures on a fully
unstructured triangulation. \Cref{fig:ffs} shows the computed density and
pressure fields at $t=4$.

Beyond this case, the distribution ships a broad battery of verification
problems, all defined in 

\texttt{src/Cases/Test\_Cases.py}: the full set of
nineteen two-dimensional Riemann configurations of \citet{lax1998solution}; the
double Mach reflection \citep{woodward1984numerical}; the Kelvin-Helmholtz
instability; the Taylor-Green vortex; several isentropic low-Mach vortex problems
(advected, dipole, co-rotating, and merging Lamb-Oseen pairs, initialised in
mechanical equilibrium via a pressure-Poisson projection to suppress spurious
acoustics); and a cylinder-in-channel geometry. Structure-preserving behaviour is
monitored through the entropy, energy, and enstrophy diagnostics noted above.

\section{Future work}
\label{sec:future_work}
One of the main limitations of the current implementation is that the mesh input
is restricted to the built-in MeshPy/Triangle\footnote{\url{https://github.com/inducer/meshpy}} generator. A VTK mesh reader is
planned, which will allow the solver to operate on arbitrary unstructured meshes
produced by external mesh generators. The current implementation is also limited to
two-dimensional triangular meshes; extension to three-dimensional tetrahedral
meshes is planned. Finally, the current implementation does not include a
turbulence model; the addition of a turbulence closure, such as a Smagorinsky
model or a Reynolds-averaged Navier-Stokes (RANS) model, is planned for future releases.

\section{Conclusion}
\label{sec:conclusion}
JAX-FVM provides, in a single compact code base, an entropy-stable compressible
Euler/Navier-Stokes finite volume solver on unstructured meshes that is
differentiable end to end and runs on CPU or GPU. Its distinctive combination of
unstructured geometry, entropy-consistent numerics, and JAX-based automatic
differentiation makes it a natural platform for gradient-based inverse problems,
PDE-constrained optimisation, and hybrid physics/machine-learning modelling--a
space that current structured-grid differentiable CFD frameworks do not cover.
Planned extensions include gradient-based shape and inverse-design workflows built
on the existing \texttt{flax}/\texttt{optax} dependencies, and completion of the
VTK mesh-input path.

\section*{Acknowledgements}
This work was carried out at INRIA. A deep thank you to Julien, Anass and Mathias for their support and feedback.

\section*{Declaration of Generative AI and AI-assisted technologies in the
writing process}
During the preparation of this work the authors used the Claude
large language model inside the VSCode IDE to assist in writing code.
After using this tool, the authors reviewed and edited the content as
needed and takes full responsibility for the content of the published
article.

\newpage
\bibliographystyle{elsarticle-harv}
\bibliography{bibliography}

\end{document}